\documentclass[a4paper,reqno,12pt]{amsart}
\usepackage{amsmath,amssymb,amscd}

\usepackage{type1cm}
\usepackage[T1]{fontenc}
\usepackage{textcomp}

\usepackage{bbm}
\usepackage{yhmath,mathrsfs}

\numberwithin{equation}{section}	

\theoremstyle{plain}
\newtheorem{thm}{Theorem}[section]	
\newtheorem*{thm*}{Theorem}

\newtheorem{lem}[thm]{Lemma}

\newtheorem{prop}[thm]{Proposition}

\newtheorem{cor}[thm]{Corollary}

\theoremstyle{definition}

\newtheorem*{ex*}{Example}

\theoremstyle{remark}
\newtheorem{rem}[thm]{Remark}

\theoremstyle{remark}

\input ncpfaff.def

\title%
[Minor summation formula of noncommutative Pfaffians]%
{A central element in the universal enveloping algebra 
 of type $\mathsf{D}_n$ \\ 
 via minor summation formula of Pfaffians}
\author{Takashi Hashimoto%
}
\address{
  Department of Information and Knowledge Sciences, 
  Faculty of Engineering, Tottori University, 
  4-101, Koyama-Minami, Tottori, 680-8552, Japan
    }
\email{thashi@ike.tottori-u.ac.jp}
\date{\today}

\begin{document}

\begin{abstract}
It is known that
the universal enveloping algebra $U(\o_{2n})$ of
the orthogonal Lie algebra $\o_{2n}$ of type $\mathsf{D}_n$
has a central element described in terms of Pfaffian
of a certain matrix which is alternating along the anti-diagonal
(we call such a matrix \textit{anti-alternating} for short)
whose entries are in $U(\o_{2n})$. 
In this paper, 
we establish minor summation formulae of Pfaffian
for the noncommutative anti-alternating matrix,
as well as for commutative anti-alternating matrix.
As an application, 
we show that 
the eigenvalues of the Pfaffian-type central element 
on the irreducible representations of $\o_{2n}$
can be easily computed.
\end{abstract}

\maketitle

\section{Introduction}
\label{intro}

Let us denote by $\o_{2n}$ the complex orthogonal Lie algebra 
defined by
\[
 \o_{2n}=\{ X \in \Mat{2n}{\C}; \tp{X} J_{2n} + J_{2n} X = O \},
\]
where $J_{2n}$ is the nondegenerate symmetric matrix 
with $1$'s on the anti-diagonal and $0$'s elsewhere.
Then a matrix $X$ is in $\o_{2n}$ if and only if 
it is alternating along the anti-diagonal,
which we call \textit{anti-alternating} for short in this paper.
We remark that 
the subspaces of $\o_{2n}$ consisting of the diagonal matrices
and of the upper triangular matrices
are Cartan subalgebra and 
nilpotent subalgebra spanned by positive root vectors,
which we denote by $\mathfrak{h}$ and $\mathfrak{n}$,
respectively.

It is known that the center $ZU(\o_{2n})$ of
the universal enveloping algebra $U(\o_{2n})$
of the Lie algebra $\o_{2n}$ is generated by elements
that are described in terms of 
column (minor) determinants and Pfaffian 
of certain matrices 
$\bs{\Phi}=(\Phi_{i,j})_{i,j=1,\dots,2n}$
whose diagonal and upper triangular entries are 
in $U(\mathfrak{h})$ and $\mathfrak{n}$, respectively
(see below for the Pfaffian type; 
one must shift diagonal entries 
for the determinant type.)
%
Here, 
for a matrix $\bs{\Phi}=(\Phi_{i,j})_{i,j}$
with $\Phi_{i,j}$ in
a noncommutative associative algebra in general,
the column determinant of $\bs{\Phi}$,
which we denote by $\det( \bs{\Phi} )$, 
is defined to be
\begin{equation}  
\label{e:nc_det}
 \det( \bs{\Phi} )
  =\sum_{\sigma \in \mathfrak{S}_{2n}}
       \sgn{\sigma} \Phi_{\sigma(1),1} \Phi_{\sigma(2),2} 
                 \cdots \Phi_{\sigma(2n),2n}.
\end{equation}
Minor column determinants are defined similarly.
If, moreover, $\bs{\Phi}$ is anti-alternating, 
by which we mean that $\bs{\Phi} J_{2n}$ is alternating as above,
the Pfaffian of $\bs{\Phi} J_{2n}$,
which we denote by $\Pf{\bs{\Phi}}$ simply, 
is defined to be
\[
 \Pf{\bs{\Phi}}
   = \frac{1}{2^n n!} \sum_{\sigma \in \mathfrak{S}_{2n}}
     \sgn{\sigma}\, 
         \tilde{\Phi}_{\sigma(1),\sigma(2)} \tilde{\Phi}_{\sigma(3),\sigma(4)}
                 \cdots   \tilde{\Phi}_{\sigma(2n-1),\sigma(2n)},
\]
where $\tilde{\Phi}_{i,j}$ denotes 
the $(i,j)$th entry of $\bs{\Phi} J_{2n}$.

It is easy to compute the eigenvalues of the central elements 
of determinant type on the irreducible representations 
of $\o_{2n}$ with highest weight $\lambda$;
if we apply the $\det(\bs{\Phi})$ to the highest vector,
the only term that survives in the sum \eqref{e:nc_det}
is the one that corresponds to $\sigma=1$
since $\Phi_{i,j}$ is in $\mathfrak{n}$ if $i<j$
as we remarked above (cf. \cite{Itoh05}).

The same principle will not work for the Pfaffian-type element.

Let us first consider the commutative case.
Write an anti-alternating matrix $X \in \o_{2n}$ as 
\begin{equation}   
\label{e:coloring_intro}
 X=\begin{bmatrix}
   a & b \\
   c & -J_n \tp{a} J_n
  \end{bmatrix},
\end{equation}
with $a,b,c$ all $n \times n$ matrices such that
$b,c$ are anti-alternating.
Then we find that
the Pfaffian of $X J_{2n}$,
which we denote by $\Pf{X}$ simply as above,
is expanded as follows 
(Corollary \ref{c:msf_square}):
\begin{equation}     
\label{e:msf_square_intro}
 \Pf{X} = \sum_{k=0}^{ \lfloor n/2 \rfloor } 
          \sum_{ \substack{
                 I, J \subset [n] \\%
                 |I|=|J|=2k}  }
           \sgn{\bar{I},I} \sgn{\bar{J},J}
           \det(a^{\bar I}_{\bar J}) \Pf{b_{I}} \Pf{c_{J}}
\end{equation}
with the second sum over all index sets $I$ and $J$,
both of cardinality $k$ and contained in $[n]:=\{1,2,\dots,n\}$,
where $\bar{I}$ and $\bar{J}$ denote the complements of 
$I$ and $J$ in $[n]$ respectively,
$a^{\bar{I}}_{\bar{J}}$ submatrix of $a$
whose row and column indices are 
in $\bar{I}$ and $\bar{J}$ respectively,
and $b_I, c_I$ submatrices of $b,c$ 
whose row and column indices are both in $I$
(see Section 3 for details.)

Furthermore,
the minor summation formula of Pfaffian \eqref{e:msf_square_intro}
holds true for a rectangular submatrix $a$.
More precisely,
for positive integers $p,q$ with $p+q=2n$,
the formula holds true 
if we write an anti-alternating matrix $X \in \o_{2n}$
as in \eqref{e:coloring_intro},
but in this case,
with $a,b,c$ of size $p \times q, p \times p, q \times q$,
respectively, 
and $(2,2)$-block replaced by $-J_q \tp{a} J_p$ 
(Theorem \ref{th:msf}).
It is immediate to see that
the minor summation formula given in \cite[Theorem 3.5]{IW06}
corresponds to ours with $p$ and $q$ both even.

Now let us return to the noncommutative case.
Let $\bs{X}=(X_{i,j})$ be an anti-alternating matrix
whose $(i,j)$th entry is given by
$X_{i,j}:=E_{i,j} -J_{2n} E_{j,i} J_{2n} %
 \in \o_{2n} \subset U(\o_{2n})$,
where $E_{i,j}$ is the matrix unit 
with $1$ in $(i,j)$th entry and 0 elsewhere.
Then we find that the following expansion formula of
the noncommutative Pfaffian $\Pf{\bs{X}}$ holds true,
which is our main result (Theorem \ref{th:main_thm})
in this paper:
\begin{thm*}
If we color $\bs{X}$ as in \eqref{e:coloring_intro},
then the Pfaffian $\Pf{\bs{X}}$ is expanded as follows:
\begin{equation*} 
  \Pf{\bs{X}}
   = \sum_{k=0}^{ \lfloor n/2 \rfloor } 
            \sum_{ \substack{I, J \subset [n] \\%
                   |I|=|J| = 2k} }
       \sgn{\bar{I},I} \sgn{\bar{J},J}
        \det \left(  a^{\bar{I}}_{\bar{J}} \,
               +  \mathbbm{1}^{\bar I}_{\bar J} \,\, \bs{\rho}(|\bar J|)
             \right) 
                  \Pf{c_{J}} \Pf{b_{I}},
\end{equation*}
where $\mathbbm{1}$ denotes the $n \times n$-identity matrix,
and $\bs{\rho}(j)=\diag(j-1,j-2,\cdots,0)$. 
\end{thm*}
Using the formula,
one can easily compute the eigenvalues of $\Pf{\bs{X}}$ 
on the irreducible representations of $\o_{2n}$
as in the case of the determinant-type elements described above,
since $b_{i,j} \in \mathfrak{n}$ for all $i,j$.

\begin{ex*}
If $n=2$, then
a matrix $X \in \o_{2n}$ 
colored as in \eqref{e:coloring_intro}
looks like
\begin{equation}
\label{e:ex_1}
 X = \begin{bmatrix}
       a_{1,1} & a_{1,2} & b_{1,2} & 0       \\
       a_{2,1} & a_{2,2} & 0       &-b_{1,2} \\
       c_{1,2} & 0       & -a_{2,2}&-a_{1,2} \\
       0       &-c_{1,2} & -a_{2,1}&-a_{1,1}
      \end{bmatrix}
\end{equation}
and the Pfaffian $\Pf{X}$ is given by
\begin{equation}
\label{e:ex_2}
 \Pf{X} =  a_{1,1}a_{2,2} -a_{2,1}a_{1,2}  + \; c_{1,2} b_{1,2},
\end{equation}
of which the first and the second terms corresponds 
to $I=J=\varnothing$ and the last one to $I=J=\{1,2\}$ 
in the right-hand side of \eqref{e:msf_square_intro}, respectively.

Passing to the noncommutative case,
let us write the matrix $\bs{X}$ as in \eqref{e:ex_1}.
Noting that
symmetrization of the right-hand side of \eqref{e:ex_2}
yields $\Pf{\bs{X}}$ by definition, 
we obtain 
\begin{align*}
 \Pf{\bs{X}} 
  &=\frac{1}{2}%
      \Bigl(
            a_{1,1} a_{2,2} - a_{2,1} a_{1,2}  +  c_{1,2} b_{1,2}
            + a_{2,2} a_{1,1} - a_{1,2} a_{2,1} +  b_{1,2} c_{1,2} 
      \Bigr)
         \\
  &= \left(a_{1,1}+1 \right)a_{2,2}-a_{2,1}a_{1,2} 
      + c_{1,2}b_{1,2},
\end{align*}
where we used the commutation relations 
\eqref{e:CR_abc} below.
\end{ex*}

The paper is organized as follows:
In Sect.2, for the sake of completeness, 
we include the proof for the fact that 
one can construct a Pfaffian-type central element 
of the universal enveloping algebra 
of the orthogonal Lie algebra 
in an \textit{arbitrary} realization.
Namely, one may take 
an arbitrary nondegenerate symmetric matrix $S$ of size $2n \times 2n$
in stead of $J_{2n}$ in the definition of $\o_{2n}$ above.
Section 3 is devoted to the proof of 
the commutative minor summation formula 
of the Pfaffian for $X \in \o_{2n}$ 
colored as in \eqref{e:coloring_intro} 
with $a$ rectangular submatrix of size $p \times q$.
We give its proof by induction on $p+q$, the size of $X$,
since it reveals that
iteration of the row/column expansion formula of the Pfaffian 
yields our formula.
In Sect.4,
we prove the theorem, i.e., 
the noncommutative minor summation formula 
of the Pfaffian $\Pf{\bs{X}}$ for the matrix $\bs{X}$,
and show that the eigenvalues of the central element $\Pf{\bs{X}}$
on the highest weight modules can be easily computed.
Finally, 
in the Appendix
we will give another proof of 
the commutative minor summation formula
using the exterior calculus.

\subsection*{Convention}
For a positive integer $n$,
let us denote by $[n]$ the set $\{ 1,2,\dots,n \}$,
$[-n]$ the set $\{ -n,\dots,-2,-1 \}$,
and $[\pm n]$ the union $[n] \cup [-n]$.
For a pair of index sets $I \subset J$,
its complement $\bar I$ is always taken in $J$, 
unless otherwise mentioned.
The symbol $\sqcup$ denotes the disjoint union.
For index sets 
\( I=\{ i_1 <\cdots <i_r \}, J=\{ j_1  <\cdots <j_s \} \)
and
\( K=\{ k_1  <\cdots <k_{r+s} \} \)
such that $K=I \sqcup J$,
let us denote by 
\(
\sgn{\begin{matrix}   K   \\
                     I,J
     \end{matrix}
    }
\)
the signature of the permutation 
\[
\begin{pmatrix} 
        k_1  & \dots & k_r  & k_{r+1} & \dots & k_{r+s}  \\
        i_1  & \dots & i_r  & j_1     & \dots & j_s
 \end{pmatrix}.
\]
%
%
When dealing with the Pfaffian  
of an anti-alternating matrix of size $2n$,
it is convenient to use the signed index. 
Namely, for any index $i \in [2n]$, 
we shall agree that $-i$ stands for $2n+1-i$.
Finally, for a real number $x$,
$\floor{x}$ denotes the greatest integer not exceeding $x$.

\section{Noncommutative Pfaffian}
\label{s:nc_Pfaffian}

Let $\bs{A}=(A_{i,j})_{i,j \in [2n]}$, $A_{j,i}=-A_{i,j}$, 
be an alternating matrix of size $2n \times 2n$ 
whose entries are elements of 
an associative $\C$-algebra $\mathscr A$.
Then the Pfaffian of $\bs{A}$, denoted by $\Pf{\bs{A}}$,
is defined to be
\begin{align}
 \Pf{\bs{A}}
  &= \frac{1}{2^n n!} \sum_{\sigma \in \mathfrak{S}_{2n}}
     \sgn{\sigma}\, A_{\sigma(1),\sigma(2)} A_{\sigma(3),\sigma(4)}
                     \cdots   A_{\sigma(2n-1),\sigma(2n)}
       \notag    \\
  &= \frac{1}{n!} \hspace{-10pt}
        \sum_{ \begin{subarray}{c}
                 \sigma \in \mathfrak{S}_{2n} \\
                 \sigma(2i-1) < \sigma(2i)
               \end{subarray}
             }  \hspace{-10pt}
     \sgn{\sigma}\, A_{\sigma(1),\sigma(2)} A_{\sigma(3),\sigma(4)}
                     \cdots   A_{\sigma(2n-1),\sigma(2n)}.
        \label{e:def_ncPfaffian}
\end{align}
If the algebra $\mathscr{A}$ is commutative,
the above definition reduces to 
\[
 \Pf{\bs{A}}
  = \sum_{\sigma}
     \sgn{\sigma}\, A_{\sigma(1),\sigma(2)} A_{\sigma(3),\sigma(4)}
                     \cdots   A_{\sigma(2n-1),\sigma(2n)},
\]
where the sum is taken over those $\sigma$ satisfying
\[ \sigma(2i-1) <\sigma(2i) \text{ for } i=1,2,\dots,n
   \text{ and }
 \sigma(1) <\sigma(3) <\cdots <\sigma(2n-1).
\]
For $g \in \GL{2n}(\C)$,
one can show that 
\begin{equation}
\label{e:equivariance_of_Pf}
 \Pf{g \bs{A} \tp{g}}=\det g \, \Pf{\bs{A}}
\end{equation} 
even if $\mathscr{A}$ is noncommutative,
as well as if $\mathscr{A}$ is commutative (see \cite{IU01}.)

Let $S$ be a nondegenerate symmetric matrix of size $2n \times 2n$.
We define the orthogonal Lie group $\O(\C^{2n},S)$ 
and its Lie algebra $\o(\C^{2n},S)$ by
\begin{align*}
 \O(\C^{2n},S) &:= \{ g \in \GL{2n}(\C); \tp{g} S g =S \},
        \\
 \o(\C^{2n},S) &:= \{ X \in \Mat{2n}{\C}; \tp{X} S + S X=O \}.
\end{align*}
Throughout this section, 
we set $G:=\O(\C^{2n},S)$ and $\g:=\o(\C^{2n},S)$ for brevity.
Then one can construct an element 
belonging to the center $ZU(\g)$  
of the universal enveloping algebra $U(\g)$ 
in terms of Pfaffian as follows.
Set $X_{i,j}:=E_{i,j}-S^{-1} \tp{E}_{i,j}S$.
Clearly, $X_{i,j} \in \g \subset U(\g)$.
Define a matrix $\bs{X}$ 
with entries in $U(\g)$ 
by 
\[
 \bs{X}:=(X_{i,j})_{i,j \in [2n]}
 = \sum_{i,j \in [2n]} E_{i,j} \otimes X_{i,j} 
     \in \Mat{2n}{\C} \otimes U(\g).
\]

\begin{lem}
The matrix $\bs{X}S$ is alternating.
Namely,
\begin{equation*}
 \bs{X} S + S \tp{\bs{X}} = O.
\end{equation*}
\end{lem}
\begin{proof}
This is a straightforward calculation, and is left to the reader.
\qed
\end{proof}
Thus one can define a Pfaffian of 
the alternating matrix $\bs{X} S$,
which we denote by $\Pf{\bs{X}}$ simply,
when there is no danger of confusion.

Following \cite{IU01},
let us write
\[
 \Ad_g(\bs{Y}) := \left( \Ad(g) Y_{i,j} \right)_{i,j \in [2n]}
     =\sum_{i,j} E_{i,j} \otimes \Ad(g)Y_{i,j}
\] 
for $g \in \GL{N}(\C)$ and 
a matrix 
$\bs{Y}=(Y_{i,j})_{i,j}=\sum_{i,j} E_{i,j} \otimes Y_{i,j}%
  \in \Mat{2n}{\C} \otimes U(\gl_N)$.
 
\begin{lem}
\label{l:Adjoint}
For $g \in G$,
we have
\[
 \Ad_g(\bs{X}) = \tp{g} \bs{X} \tp{g}^{-1}.
\]
\end{lem}
\begin{proof}
Setting $g=(g_{a b})$ and $g^{-1}=(g^{a b})$,
one has
\begin{align*}
\Ad(g) E_{i,j} 
  &= \sum_{a,b,c,d}  
        g_{a b} g^{c d} 
              E_{a,b} E_{i,j} E_{c,d}
   = \sum_{a, b} g_{a i} g^{j b} E_{a,b}.
\end{align*}
Similarly,
setting $S=(s_{a b})$ and $S^{-1}=(s^{a b})$,
one has 
\begin{align*}
\Ad(g) \left( S^{-1} \tp{E_{i,j}} S \right)
 &= \sum_{a,b,k,l} s_{i l} g^{l b} g_{a k} s^{k j} E_{a,b}.
\end{align*}
Therefore, 
denoting by $\bs{E}$ 
the matrix whose $(i,j)$th entry is $E_{i,j}$,
one has
\begin{align*}
\Ad_g( \bs{X} )
 &= \tp{g} \bs{E} \tp{g}^{-1} - S g^{-1} \tp{\bs{E}} g S^{-1} 
        \\
 &= \tp{g} \left( \bs{E} - S \tp{\bs{E}} S^{-1} \right) \tp{g}^{-1}
        \\
 &= \tp{g} \bs{X} \tp{g}^{-1}
\end{align*}
since 
$\left( S^{-1} \tp{E_{i,j}} S \right)_{i,j} = S \tp{\bs{E}} S^{-1}$
and $g \in G$.
\qed
\end{proof}

%
%
\begin{prop}
\label{p:Pf_is_central}
$\Pf{\bs{X}}$ belongs to the center $ZU(\g)$
of the universal enveloping algebra $U(\g)$.
More precisely, it satisfies 
\begin{equation}
 g \Pf{\bs{X}} g^{-1}=\det g \Pf{\bs{X}}
\end{equation}
for all $g \in G$.
\end{prop}
\begin{proof}
First we note that $\Ad_g(\bs{X} S)= \Ad_g(\bs{X}) S$.
In fact,
\begin{align*}
 \Ad_g(\bs{X}S) 
  &= \Ad_g \left( 
        \sum E_{i,j} \otimes X_{i,j} \cdot \sum E_{k,l} \otimes s_{k l}
            \right)
     \\
  &= \sum E_{i,j} \otimes \Ad(g)X_{i,j} 
                      \cdot \sum E_{k,l} \otimes s_{k l}
     \\
  &= \Ad_g(\bs{X}) S, 
\end{align*}  
whence $\Ad_g(\bs{X}S)=\tp{g} \bs{X} \tp{g}^{-1} S$ 
by Lemma \ref{l:Adjoint}.
Therefore,
setting $\bs{X}S =: \sum E_{i,j}\otimes \tilde{X}_{i,j}$,
it follows from \eqref{e:equivariance_of_Pf} that
\begin{align*}
 g \Pf{\bs{X}} g^{-1}
  &= \frac{1}{2^n n!}  
            \sum_{\sigma \in \mathfrak{S}_{2n}}
              \Ad(g) \tilde{X}_{\sigma(1),\sigma(2)} 
                   \cdots \Ad(g) \tilde{X}_{\sigma(2n-1),\sigma(2n)}
        \\
  &= \Pf{ \tp{g} \bs{X} \tp{g}^{-1} S}
     = \Pf{ \tp{g} \bs{X} S g}
        \\
  &= \det g \Pf{\bs{X}}, 
\end{align*}
since $g \in G$.
\qed
\end{proof}

\section{Commutative Minor Summation Formula}

For a positive integer $N$,
let us denote by $J_{N}$ the $N \times N$ matrix 
with $1$'s on the anti-diagonal and $0$'s elsewhere:
\[
 J_{N}:= \begin{bmatrix}  
            &        & 1 \\[-4pt]
            & \adots &   \\[-4pt]
          1 &        &
         \end{bmatrix}.
\] 
Take $J_{2n}$ as the nondegenerate symmetric matrix $S$ 
in Section \ref{s:nc_Pfaffian},
and denote the orthogonal group $\O(\C^{2n},J_{2n})$ 
and its Lie algebra $\o(\C^{2n},S)$ by $\O_{2n}$ and $\o_{2n}$,
respectively:
\begin{align*}
\O_{2n} &:= \{ g \in \GL{2n}(\C); \tp{g} J_{2n} g=J_{2n} \}, 
     \\
\o_{2n} & :=\{ X \in \Mat{2n}{\C}; \tp{X} J_{2n} + J_{2n} X =0 \}.
\end{align*}
Note that a matrix $X$ is in $\o_{2n}$
if and only if it is alternating along the anti-diagonal,
which we call \textit{anti-alternating} for short,
so that one can define the Pfaffian of $XJ$,
which we denote by $\Pf{X}$ simply,
as in the previous sections.

%
%
\begin{rem}
In general, 
for a symmetric matrix $S$,
$X$ is in $\o(\C^{2n},S)$ if and only if 
$X S^{-1}$ is alternating.
Note that $J_{2n}=J_{2n}^{-1}$ 
in our case of anti-alternating matrices.
\end{rem}

For positive integers $p,q $ with $p+q=2n$,
we write a matrix $X \in \o_{2n}$ as
\begin{equation}
\label{e:coloring}
  X= \begin{bmatrix}
        a &   b             \\ 
        c & -J_q \tp{a} J_p
    \end{bmatrix},
\end{equation}
where $a,b,c$ are matrices of size 
$p \times q,p \times p,q \times q$ respectively, 
such that $b,c$ are anti-alternating, 
and then parameterize $a,b,c$ as 
\begin{equation}
\label{e:parameterize}
a=\sum_{i \in [p], j \in [q]} a_{i,j} E_{i,j},
   \quad
b=\sum_{i,j \in [p]} b_{i,j} E_{i,-j},
  \quad
c=\sum_{i,j \in [q]} c_{i,j} E_{-j,i}
\end{equation}
where $a_{i,j},b_{i,j},c_{i,j} \in \C$ 
such that $b_{j,i}=-b_{i,j}$ and $c_{j,i}=-c_{i,j}$.
Here $E_{i,j}$ denotes the matrix unit
taking the basis vector $e_{j}$ to $e_{i}$,
where $\{ e_i \}_{i \in [\pm n]}$ is 
the standard basis for $\C^{2n}$.
Now we define their submatrices by 
\[
 a^{I}_{J}:=(a_{i,j})_{i \in I,j \in J},
   \quad
 b_{I}:=(b_{i,j})_{i,j \in I},
   \quad
 c_{J}:=(c_{i,j})_{i,j \in J}
\]
for $I \subset [p], J \subset [q]$.
Note that $b_I$ and $c_J$ are still anti-alternating.

\begin{thm}
\label{th:msf}
Let $r=\min(p,q)$ and $\epsilon$  the parity of $p$,
i.e., is equal to $0$ or $1$ 
according as $p$ is even or odd.
If
\(
X=\begin{bmatrix} 
 a & b \\ 
 c & -J_q \tp{a} J_p
 \end{bmatrix}
\)
is a matrix in $\mathfrak{o}_{2n}$ 
with $a,b,c$ parameterized as in \eqref{e:parameterize},
then the Pfaffian $\Pf{X}$ is expanded as follows:
\begin{equation}
\label{e:msf}
   \Pf{X} = \sum_{k=0}^{ \lfloor r/2 \rfloor } 
            \sum_{ \substack{I \subset [p],J \subset [q] \\%
                   |\bar I|=|\bar J| = 2 k + \epsilon} }
           \sgn{\bar{I},I} \sgn{\bar{J},J}
           \det(a^{\bar I}_{\bar J}) \Pf{b_{I}} \Pf{c_{J}},
\end{equation}
where 
\(
\sgn{\bar I,I}
 =\sgn{ \begin{matrix} 
          1   ... p  \\
          \bar I, I
        \end{matrix}
      } 
\)
and 
\(
\sgn{\bar J,J}
 =\sgn{ \begin{matrix} 
          1  ... q  \\
          \bar J, J
        \end{matrix}
      }. 
\) 
\end{thm}

%
%
\begin{rem}
For an alternating matrix $A=(a_{i,j})_{i,j \in [2n]}$
and an indices $I \subset [2n]$,
let us denote by $A_I$ the alternating submatrix 
consisting of $(a_{i,j})_{i,j \in I}$, i.e., 
of entries whose row and column indices are both in $I$.
For $i,j \in [2n]$,
the $(i,j)$th cofactor Pfaffian $\gamma_{i,j}(A)$ 
is then defined to be 
\[
 \gamma_{i,j}(A)
  =\begin{cases}
    (-1)^{i+j-1} \Pf{A_{ [1..\hat{i}..\hat{j}..2n] }}
        &\text{if $i<j$},   \\
    \hspace{14pt} 0
        &\text{if $i=j$},  \\
    (-1)^{i+j} \Pf{A_{ [1..\hat{j}..\hat{i}..2n] }}
        &\text{if $i>j$},   \\
   \end{cases}
\]
where $\hat{i}$ means omitting $i$.
Then, as in the case of determinant,
the following expansion formula holds:
\begin{equation}
\label{e:coPfaffian_expansion}
 \delta_{i,j} \Pf{A} = \sum_{k=1}^{2n} a_{i,k} \, \gamma_{j,k}(A)
\end{equation}
(see \cite{Hirota92}.)

The co-Pfaffian matrix of $A$ is,
by definition,
an alternating matrix 
whose $(i,j)$th entries are given by
$ \gamma_{i,j}(A)$ for $i,j \in [2n]$,
which we denote by $\hat{A}$.
Then one can verify that
\begin{equation}  
\label{e:relation_with_IW06}
 \frac{ \Pf{A_I} }{ \Pf{A} }
   =\sgn{I,\bar{I}}  \Pf{({\hat{A}}/{\Pf{A}})_{\bar I}}
\end{equation}
if $A$ is invertible.
Note that \eqref{e:relation_with_IW06} makes sense 
only if $|I|$ is even.
Using the relation \eqref{e:relation_with_IW06}, 
it is immediate to see that the minor summation formula 
given in \cite[Theorem 3.5]{IW06} 
coincides with ours \eqref{e:msf} with 
$p$ and $q$ both even.
\end{rem}

%
%
%
%
\begin{proof}[Proof of Theorem \ref{th:msf}] 
Here we will give the proof using induction on $p+q$,
since it reveals that 
iteration of the co-Pfaffian expansion 
of the Pfaffian \eqref{e:coPfaffian_expansion}
yields our minor summation formula.

It suffices to prove when $p \leqsl q$.
The case where $p=0,q=1$ is trivial. 
Now expanding $\Pf{X}$ of the matrix $X$ given by 
\eqref{e:coloring} and \eqref{e:parameterize} 
with respect to the first row,
one obtains
\begin{equation}
\label{e:cofactor_expansion}
 \Pf{X} = \sum_{j=1}^{q} a_{1,j} \gamma_{1,-j} 
           +\sum_{j=2}^{p} b_{1,j} \gamma_{1,j}, 
\end{equation}
where we put $\gamma_{i,j}=\gamma_{i,j}(X J_{2n})$ for brevity.
Note that $\gamma_{i,j}$ is given by, up to the sign,
the Pfaffian of the submatrix obtained by deleting 
the $i$th and $j$th rows, 
and $-i$th and $-j$th columns from $X$,
which is anti-alternating.

By inductive hypothesis,
one obtains that
\begin{align}
 \gamma_{i,-j} 
   &= (-)^{1+j} \sum_{k=0}^{r-1}
                 \sum_{
                    \begin{subarray}{c}
                      I \subset [2...p] \\
                      J \subset [1..\hat{j}..q] \\
                      |\bar I|=|\bar J|= 2k +\epsilon_1
                    \end{subarray}
                  }
       \sgn{ \begin{matrix}
               2...p      \\  
               I, \bar{I}
             \end{matrix} 
         }
       \sgn{ \begin{matrix}
               1..\hat{j}..q \\
               J, \bar J  
            \end{matrix}
         }
       \det ( a^{\bar{I}}_{\bar{J}} ) \Pf{b_{I}} \Pf{c_{J}},
     \label{e:gamma_{i,-j}} \\
 \gamma_{i,j} 
   &= (-)^{j} \sum_{k=0}^{r-1}
                 \sum_{
                    \begin{subarray}{c}
                      I \subset [2..\hat{j}..p] \\
                      J \subset [q] \\
                      |\bar I|=|\bar J|= 2k +\epsilon
                    \end{subarray}
                  }
       \sgn{ \begin{matrix}
               2..\hat{j}..p     \\  
               I, \bar I
            \end{matrix}
         }
       \sgn{ \begin{matrix}
             1...q \\
             J, \bar J  
            \end{matrix}
         }
       \det ( a^{\bar{I}}_{\bar{J}} ) \Pf{b_{I}} \Pf{c_{J}},
     \label{e:gamma_{i,j}}
\end{align}
where 
$\epsilon_1$ is the parity of $p-1$,
i.e., $\epsilon_1=1-\epsilon$.

In the right-hand side of \eqref{e:cofactor_expansion},
we denote 
the term of degree $2k+\epsilon$ in the variables $a_{i,j}$
by $T_k$ ($k=0,1,\dots,r$).
Then it follows from \eqref{e:gamma_{i,-j}} 
and \eqref{e:gamma_{i,j}} that 
\begin{align}
\label{e:cofactor_expansion2}
 T_k &= \sum_{j=1}^{q} (-)^{1+j} 
                 \sum_{
                    \begin{subarray}{c}
                      I \subset [2...p] \\
                      J \subset [1..\hat{j}..q] \\
                      |\bar I|=|\bar J|= 2k-\epsilon_1
                    \end{subarray}
                  }
       \sgn{ \begin{matrix}
               2...p      \\  
               I, \bar{I}
             \end{matrix} 
         }
       \sgn{ \begin{matrix}
               1..\hat{j}..q \\
               J, \bar J  
             \end{matrix}
         }
       a_{1,j} \det ( a^{\bar{I}}_{\bar{J}} ) \Pf{b_{I}} \Pf{c_{J}}
         \notag   \\ 
   &\phantom{=} +\sum_{j=2}^{p} (-)^{j}
                 \sum_{
                    \begin{subarray}{c}
                      I \subset [2..\hat{j}..p] \\
                      J \subset [q] \\
                      |\bar I|=|\bar J|= 2k+\epsilon
                    \end{subarray}
                  }
       \sgn{ \begin{matrix}
              2..\hat{j}..p     \\  
              I, \bar I
             \end{matrix}
         }
       \sgn{ \begin{matrix}
              1...q \\
              J, \bar J  
             \end{matrix}
         }
       \det ( a^{\bar{I}}_{\bar{J}} ) \,b_{1,j} \Pf{b_{I}} \Pf{c_{J}}.
\end{align}
Let us rewrite the first sum in the right-hand side 
of \eqref{e:cofactor_expansion2} as 
\[
 \sum_{\begin{subarray}{c}
         I \subset [2...p] \\
         |\bar I|=2k-\epsilon_1
       \end{subarray}
     }
    \sgn{\begin{matrix}
           2...q \\
           I, \bar I 
         \end{matrix}
      }    \Pf{b_{I}}
     \Biggl(
       \sum_{j=1}^{q}
       \sum_{\begin{subarray}{c}
              J \subset [1..\hat{j}..q] \\
              |\bar J|=2k-\epsilon_1
             \end{subarray}
           }
        (-)^{j+1} \sgn{\begin{matrix}
                         1..\hat{j}..q \\
                         J, \bar J 
                       \end{matrix}
                       }
             a_{1,j} \det ( a^{\bar{I}}_{\bar{J}} ) \Pf{c_{J}}
     \Biggr).
\]
In the brace of the equation above,
fixing an index set $J_1 \subset [q]$ 
of length $q-1-(2k-\epsilon_1)=q-2k-\epsilon$
and denoting its complement by ${\bar J}_1$,
one sees that the coefficient of $\Pf{c_{J_1}}$ equals
\begin{align*}
  &\sum_{j \in {\bar J}_1} (-)^{j+1} 
    \sgn{\begin{matrix}
          1..\hat{j}..q \\
          {\bar J}_1\setminus \{j\}, J_1
         \end{matrix}
       }
        a_{1,j} \det( a^{\bar I}_{ {\bar J}_1 \setminus \{j\} } ) 
           \\
 =& \sgn{\begin{matrix}
          1...q \\
          {\bar J}_1 \, J_1
         \end{matrix}}
      \det( a^{{\bar I}_1}_{{\bar J}_1} )
\end{align*}
by the expansion formula of the determinant,
where ${\bar I}_1:=\{1\} \sqcup {\bar I}$.
Hence, 
one obtains that the first sum of $T_k$ equals
\begin{equation}
\label{e:1st_sum_of_T_k}
 \sum_{\begin{subarray}{c}
         1 \notin I_1 \subset [p] \\
         |\bar{I}_1|= 2k+\epsilon
       \end{subarray}
     }
 \sum_{\begin{subarray}{c}
         J_1 \subset [q] \\
         |\bar{J}_1|= 2k+\epsilon
       \end{subarray}
     }
       \sgn{ \begin{matrix}
               I_1, \bar{I}_1
             \end{matrix} 
         }
       \sgn{ \begin{matrix}
             J_1, \bar{J}_1  
            \end{matrix}
         }
       \det ( a^{\bar{I}_1}_{\bar{J}_1} ) \Pf{b_{I_1}} \Pf{c_{J_1}}.
\end{equation}

Now let us turn to the second sum of $T_k$.
It can be written as
\[
 \sum_{\begin{subarray}{c}
         J \subset [q] \\
         |\bar J|=2k+\epsilon
       \end{subarray}
     }
    \sgn{\begin{matrix}
           1...q \\
           J, \bar J 
         \end{matrix}
      }    \Pf{c_{J}}
     \Biggl(
       \sum_{j=2}^{p}
       \sum_{\begin{subarray}{c}
              I \subset [2..\hat{j}..p] \\
              |\bar I|=2k+\epsilon
             \end{subarray}
           }
        (-)^{j} \sgn{\begin{matrix}
                         2..\hat{j}..p \\
                         I, \bar I 
                       \end{matrix}
                       }
             \det ( a^{\bar{I}}_{\bar{J}} ) b_{1,j} \Pf{b_{I}}
     \Biggr).
\]
In the brace of the equation above, 
similarly to the case of the first sum, 
fixing an index set ${\bar I}_1 \subset [2...p]$ 
of length $p-2-2k-\epsilon$,
and denoting its complement in $[p]$ by $I_1$,
one sees that the coefficient of $\det(a^{\bar{I}_1}_{J})$ equals
\[
 \sgn{\begin{matrix}
       1...p \\
       I_1,\bar{I}_1
      \end{matrix}
    } \Pf{b_{I_1}},
\]
whence, one obtains that 
the second sum of $T_k$ equals
\begin{equation}
\label{e:2nd_sum_of_T_k}
 \sum_{\begin{subarray}{c}
         J \subset [q] \\
         |\bar{J}|= 2k+\epsilon
       \end{subarray}
     }
 \sum_{\begin{subarray}{c}
         1 \in I_1 \subset [p] \\
         |\bar{I}_1|= 2k+\epsilon
       \end{subarray}
     }
       \sgn{ \begin{matrix}
               I_1, \bar{I}_1
             \end{matrix} 
         }
       \sgn{ \begin{matrix}
             J, \bar{J}  
            \end{matrix}
         }
       \det ( a^{\bar{I}_1}_{\bar{J}} ) \Pf{b_{I_1}} \Pf{c_{J}}
\end{equation}
by the expansion formula of the Pfaffian.

It follows from \eqref{e:1st_sum_of_T_k} and \eqref{e:2nd_sum_of_T_k}
that $T_k$ equals
\begin{align*}
 &\Biggl( 
     \sum_{\begin{subarray}{c}
            1 \notin I_1 \subset [p] \\
            |\bar{I}_1|= 2k+\epsilon
           \end{subarray}
       }
     \sum_{\begin{subarray}{c}
            J_1 \subset [q] \\
            |\bar{J}_1|= 2k+\epsilon
           \end{subarray}
     }
   + \sum_{\begin{subarray}{c}
               J_1 \subset [q] \\
               |\bar{J}_1|= 2k+\epsilon
              \end{subarray}
         }
     \sum_{\begin{subarray}{c}
            1 \in I_1 \subset [p] \\
            |\bar{I}_1|= 2k+\epsilon
           \end{subarray}
     }
 \Biggr)
       \sgn{ \begin{matrix}
               I_1, \bar{I}_1
             \end{matrix} 
         }
       \sgn{ \begin{matrix}
             J_1, \bar{J}_1  
            \end{matrix}
         }
       \det ( a^{\bar{I}_1}_{\bar{J}_1} ) \Pf{b_{I_1}} \Pf{c_{J_1}}
     \\
 =& \sum_{\begin{subarray}{c}
               I_1 \subset [p] \\
               |\bar{I}_1|= 2k+\epsilon
              \end{subarray}
         }
     \sum_{\begin{subarray}{c}
            J_1 \subset [q] \\
            |\bar{J}_1|= 2k+\epsilon
           \end{subarray}
         }
       \sgn{ \begin{matrix}
               I_1, \bar{I}_1
             \end{matrix} 
         }
       \sgn{ \begin{matrix}
             J_1, \bar{J}_1  
            \end{matrix}
         }
       \det ( a^{\bar{I}_1}_{\bar{J}_1} ) \Pf{b_{I_1}} \Pf{c_{J_1}}.
\end{align*}
This completes the proof.
\qed
\end{proof}

We separately state the case where $p=q=n$ 
of Theorem \ref{th:msf}, 
of which noncommutative version,
i.e., the case where the matrix entries are elements of
the universal enveloping algebra $U(\o_{2n})$,
will be considered in Section \ref{s:nc_MSF}:
\begin{cor}
\label{c:msf_square}
Let 
\(
X=\begin{bmatrix} 
 a & b \\ 
 c & -J_n \tp{a} J_n
 \end{bmatrix}
   \in \o_{2n}
\)
be as in the theorem,
where $a,b,c$ are all of size $n \times n$.
Then the Pfaffian $\Pf{X}$ is expanded as follows:
\begin{equation}
\label{e:msf_square}
  \Pf{X} = \sum_{k=0}^{ \lfloor n/2 \rfloor } 
          \sum_{ \substack{
                 I, J \subset [n] \\%
                 |I|=|J|=2k}  }
           \sgn{\bar{I},I} \sgn{\bar{J},J}
           \det(a^{\bar I}_{\bar J}) \Pf{b_{I}} \Pf{c_{J}}.
\end{equation}
\end{cor}

\section{Noncommutative Minor Summation Formula}
\label{s:nc_MSF}

Now we turn to the noncommutative case.

Let $\bs{A}=(A_{i,j})$ be 
an alternating matrix of size $2n \times 2n$
whose entries $A_{i,j}$ are 
in a noncommutative associative algebra $\mathscr{A}$.
The Pfaffian $\Pf{\bs{A}}$ is defined 
by \eqref{e:def_ncPfaffian}.
Note that for the alternating matrix $\bs{A}$,
if we define a $2$-form $\Theta_{\bs{A}}$ 
with coefficients in $\mathscr{A}$ by
\begin{equation*}
\Theta_{\bs{A}}
 = \sum_{i,j \in [2n]} e_{i}e_{j} A_{i,j} 
      \in \bigwedge^2 \C^{2n} \otimes \mathscr{A},
\end{equation*}
then the Pfaffian $\Pf{\bs{A}}$ is given by 
the coefficient of $e_1 e_2 \cdots e_{2n}$ 
in ${ \Theta_{\bs{A}} }^n$ divided by $2^n n!$
(\cite[Proposition 1.1]{IU01}):
\begin{equation} 
 { \Theta_{\bs{A}} }^n 
   = e_1 e_2 \cdots  e_{2n} \, 2^n n! \Pf{\bs{A}}.
\end{equation}

Let us consider the following $2n \times 2n$-matrix
whose entries are in the universal enveloping algebra $U(\o_{2n})$:
\begin{equation}
\label{e:mat_noncommutative}
  \bs{X} = (X_{i,j}) 
   \quad \text{with}  \quad
  X_{i,j} := E_{i,j}-E_{-j,-i} 
        \in \mathfrak{o}_{2n} \subset U(\mathfrak{o}_{2n})
\end{equation}
for $i,j \in [\pm n]$.
The commutation relations among $X_{i,j}$'s are given by
\begin{equation}
\label{e:commutation_relation}
 [ X_{i,j}, X_{k,l} ] 
  = \delta_{j,k} X_{i,l} + \delta_{i,l} X_{-j,-k}
   -\delta_{j,-l} X_{i,-k} - \delta_{i,-k} X_{-j,l} 
\end{equation}
for $i,j,k,l \in [\pm n]$.

Since $\bs{X}$ is anti-alternating by definition,
one can define the Pfaffian of $\bs{X} J_{2n}$,
which we also denote by $\Pf{\bs{X}}$ simply,
as usual.
Then it belongs to the center $ZU(\o_{2n})$ 
of the universal enveloping algebra,
since it is invariant under the adjoint action of
$\SO_{2n}$ by Proposition \ref{p:Pf_is_central}.

For the matrix $\bs{X}$ in \eqref{e:mat_noncommutative},
we introduce a 2-form with coefficient 
in the universal enveloping algebra $U(\o_{2n})$:
%
%
\begin{equation}
\label{e:Omega}
\Omega=\sum_{i,j \in [\pm n]} e_{i}e_{-j} X_{i,j} 
           \in  \bigwedge^2  \C^{2n} \otimes U(\o_{2n}).
\end{equation}
%
%
\begin{lem}
\label{l:Pfaffian_via_2-form}
The relation between the Pfaffian $\Pf{\bs{X}}$ and $\Omega$
is given by
\begin{equation}
\label{e:Omega2the_n}
 \Omega^n=e_1 \cdots e_n e_{-n} \cdots e_{-1} 2^n n! \Pf{\bs{X}}.
\end{equation}
\end{lem}
\begin{proof}
If we set $\tilde{\bs{X}}:=\bs{X} J_{2n}$,
its $(i,j)$th entry $\tilde{X}_{i,j}$ is given by 
$X_{i,-j}$ for $i,j \in [2n]$.
Hence
\begin{align*}
 \Omega^n 
  &= \sum_{i_1,j_1,\dots,i_n j_n \in [\pm n]}
        e_{i_1} e_{j_1} \cdots e_{i_n} e_{j_n} 
        X_{i_1,-j_1} \cdots X_{i_n,-j_n}  
               \\
  &= \sum_{i_1,j_1,\dots,i_n j_n \in [2n]}
        e_{i_1} e_{j_1} \cdots e_{i_n} e_{j_n} 
        \tilde{X}_{i_1,j_1} \cdots \tilde{X}_{i_n,j_n}
               \\
  &= e_{1} e_{2} \cdots e_{2n-1} e_{2n} 
      \sum_{\sigma \in \mathfrak{S}_{2n}}
      \sgn{\sigma} \tilde{X}_{\sigma(1),\sigma(2)} 
                      \cdots \tilde{X}_{\sigma(2n-1),\sigma(2n)}
               \\
  &= e_{1} e_{2} \cdots e_{n} e_{n+1} \cdots e_{2n-1} e_{2n} \,
      2^n n! \Pf{ \tilde{\bs{X}} } 
               \\
  &= e_{1} e_{2} \cdots e_{n} e_{-n} \cdots e_{-2} e_{-1} \,
      2^n n! \Pf{\bs{X}}, 
\end{align*}
where we use the convention that $-j=2n+1-j$ for $j \in [2n]$.
\qed
\end{proof}

Similarly to the commutative case 
with $p=q=n$ in \eqref{e:coloring},
we color the matrix $\bs{X}$ as 
\begin{equation*}   
\bs{X} =\begin{bmatrix}
          a & b \\
          c & -J_n \tp{a} J_n
        \end{bmatrix},
\end{equation*}
where $a,b,c$ are matrices of size $n \times n$ 
given by 
\begin{equation}
\label{e:nc_coloring}
a=\sum_{i, j \in [n]}  E_{i,j} \otimes a_{i,j},
   \quad
b=\sum_{i,j \in [n]}  E_{i,-j} \otimes b_{i,j},
  \quad
c=\sum_{i,j \in [n]}  E_{-j,i} \otimes c_{i,j}
\end{equation}
with 
\begin{equation}
\label{e:nc_parameterize}
 a_{i,j}=X_{i,j}, \quad
 b_{i,j}=X_{i,-j}, \quad
 c_{i,j}=X_{-j,i}.
\end{equation}
%
%
By \eqref{e:commutation_relation} and \eqref{e:nc_parameterize},
the commutation relations among $a_{i,j},b_{i,j},c_{i,j}$
are given by
\begin{align}
\label{e:CR_abc}
[a_{i,j}, a_{k,l}]
   &=  \delta_{j,k} a_{i,l} - \delta_{l,j} a_{k,j},
        \notag \\
[a_{i,j}, b_{k,l}]
   &=  \delta_{j,k} b_{i,l} - \delta_{j,l} b_{i,k},
        \notag \\
[a_{i,j}, c_{k,l}]
   &=   \delta_{i,k} c_{l,j} - \delta_{i,l} c_{k,j},
        \\
[b_{i,j}, c_{k,l}]
   &=   \delta_{j,l} a_{i,k} + \delta_{i,k} a_{j,l}
                      -\delta_{i,l} a_{j,k} - \delta_{j,k} a_{i,l},
        \notag \\
[b_{i,j}, b_{k,l}] &= [c_{i,j},c_{k,l}]  = 0
        \notag
\end{align}
for $i,j,k,l \in [n]$.

%
%
\begin{rem}
\label{r:maximal_parabolic_subalgebra}
Set $\mathfrak{h}=\bigoplus_{i \in [n]} \C a_{i,i}$.
Then $\mathfrak{h}$ is a Cartan subalgebra of $\o_{2n}$,
and each $b_{i,j}$ is a positive root vector 
with root $\epsilon_{i}+\epsilon_{j}$,
where $\epsilon_i \in \mathfrak{h}^{*}$ 
is the linear functional on $\mathfrak{h}$ 
that takes $\diag(h_1,\dots,h_n,-h_n,\dots,-h_1)$ to $h_{i}$.
The other positive root vectors are $a_{i,j}$
with root $\epsilon_{i}-\epsilon_{j}$ for $i<j$.
In fact,
setting $\mathfrak{l}:=\bigoplus_{i,j \in [n]} \C a_{i,j}$,
$\mathfrak{u}:=\bigoplus_{i,j \in [n],i<j} \C b_{i,j}$,
and  
$\mathfrak{q}:=\mathfrak{l} \oplus \mathfrak{u}$,
then $\mathfrak{q}$ is the maximal parabolic subalgebra
corresponding to the right-end root 
in Dynkin diagram of type $\mathsf{D}_n$,
with $\mathfrak{l}$ the Levi factor
which is isomorphic to $\mathfrak{gl}_n$,
and $\mathfrak{u}$ the nilradical which is abelian.
Note that $\mathfrak{u}^{-}:=\bigoplus_{i,j \in [n], i<j} \C c_{i,j}$
is also abelian
(see \cite{GW98,Knapp02}.)
\end{rem}

Corresponding to the coloring of the matrix $\bs{X}$,
we set
\begin{equation*}
\Xi = \sum_{i, j \in [n]}    e_i e_{-j} \, a_{i,j},  
     \quad  
\Theta = \sum_{i,j \in [n]}  e_i e_{j} \, b_{i,j},  
     \quad
\Theta' = \sum_{i,j \in [n]} e_{-j} e_{-i} \, c_{i,j}.  
\end{equation*}
Then obviously,
\[
 \Omega = \Theta' + 2 \Xi + \Theta.
\]

%
%
\begin{lem}  
\label{l:sl_2_triplet}
The following commutation relations hold:
\[
 [\Theta, \Theta']  = 4 \tau \Xi,  \quad 
 [\Theta, \Xi]  = 2 \tau \Theta,   \quad
 [\Theta', \Xi]  = -2 \tau \Theta',
\]
where $\tau = \sum_{i \in [n]} e_i e_{-i}$.
\end{lem}
\begin{proof}
It is straightforward to show that these relations 
follow from \eqref{e:CR_abc}.
For example, one sees that
\begin{align*}
 [\Xi,\Theta] 
  &
   = \sum_{i,j,k,l} e_i e_{-j} e_k e_l [a_{i,j}, b_{k,l}]
       \\
  &= \sum_{i,j,k,l} e_i e_{-j} e_k e_l 
            (\delta_{j,k} b_{i,l} - \delta_{j,l} b_{i,k})
       \\
  &= \sum_{i,j,l}e_i e_{-j} e_j e_l b_{i,l} 
       -  \sum_{i,j,k} e_i e_{-j} e_k e_j b_{i,k}
       \\
  &= - \sum_{i,j,l} e_{-j} e_j e_i e_l b_{i,l} 
       -  \sum_{i,j,k} e_{-j} e_j e_i e_k b_{i,k}
       \\
  &= -2\tau \Theta.
\end{align*}
The other relations follow similarly.
\qed
\end{proof}

%
%
For a parameter $u \in \C$
and $r=0,1,2,\dots$,
set 
\begin{equation}
\label{e:def_Xi(u)}
 \Xi(u) := \Xi+u \tau 
   \quad \text{and} \quad
 \Xi^{(r)}(u) := \Xi(u) \Xi(u-1) \cdots \Xi(u-r+1).
\end{equation}
The following propositions are due to 
\cite[Lemma 4.5 and Proposition 2.6]{IU01}.
%
%
\begin{prop}
\label{p:trinomial}
For $m=0,1,\dots,n$,
we have 
\begin{equation*}
\label{e:trinomial}
 \Omega^m = \sum_{ \begin{subarray}{c}
                   p,q,r \geqsl 0 \\
                   p+q+r=m
                 \end{subarray}
               }
           \frac{m!}{p! q! r!} 2^{r} \, 
             {\Xi}^{(r)}(q-p+r-1) \Theta'^{p} \Theta^{q}.
\end{equation*}
\end{prop}
\begin{proof}
By Lemma \ref{l:sl_2_triplet},
our 2-forms $\Xi,\Theta,\Theta'$ satisfy 
the same commutation relations as 
those given in \cite[Lemma 4.1]{IU01}.
Therefore, exactly the same argument therein
implies the proposition
(see \cite[Lemma 4.5]{IU01} for details.)
\qed
\end{proof}

For $j \in [n]$ and $u \in \C$,
set
%
%
\begin{equation}    
\label{e:def_eta}
 \eta_j(u)  = \sum_{i \in [n]} e_i  a_{i,j}(u)
\end{equation}
with $a_{i,j}(u)=a_{i,j} +  u \delta_{i,j}$.
Then they are anti-commutative 
\textit{when the parameter shift taken into account},
i.e., they satisfy
\begin{equation}
\label{e:anti_commutation_eta}
 \eta_{i}(u+1) \eta_{j}(u)+\eta_{j}(u+1) \eta_{i}(u)=0
\end{equation} 
for $i,j \in [n]$ (cf. \cite[Lemma 2.1]{IU01}).
Obviously, we have 
\begin{equation}
\label{e:Xi_and_eta}
 \Xi(u)= \sum_{j \in [n]} \eta_j(u) e_{-j}.
\end{equation}

Given $I,J \subset [n]$,
define submatrices of $a,b,c$ by
\[
 a^{I}_{J}:=(a_{i,j})_{i \in I,j \in J},
   \quad
 b_{I}:=(b_{i,j})_{i,j \in I},
   \quad
 c_{J}:=(c_{i,j})_{i,j \in J}
\]
as in the commutative case.
Furthermore, for economy,
we will use notations
\[
 e_I:=e_{i_1} e_{i_2} \cdots e_{i_r}
    \quad \text{and}  \quad
 e_{-J}:=e_{-j_s} e_{-j_{s-1}} \cdots e_{-j_1}
\]
if $I=\{ i_1 <i_2 <\cdots <i_r\}$ 
and $J=\{ j_1 <j_2 <\cdots j_s \}$
in what follows.

%
%
\begin{prop}    
\label{p:Xi2the_power}
For $r=0,1,\dots,n$ and $u \in \C$,
we have
\[
 \Xi^{(r)}(u+r-1)
= r! \sum_{ \begin{subarray}{c}
             I,J \subset [n] \\
             |I|=|J|=r
            \end{subarray}
        }
     e_{I} e_{-J}
     \det \left(  a^{I}_{J} \,
               +  \mathbbm{1}^{I}_{J} \diag(u+r-1,u+r-2,\dots,u)
           \right), 
\]
where 
$\det$ denotes the column determinant
and $\mathbbm{1}$ the identity matrix of size $n \times n$.
\end{prop}
\begin{proof}
First, note that the column determinant in the sum 
is explicitly given by
\begin{align}  
=& \det \begin{bmatrix}
         a_{i_1,j_1}(u+r-1) & a_{i_1,j_2}(u+r-2) & \cdots & a_{i_1,j_r}(u) \\
         a_{i_2,j_1}(u+r-1) & a_{i_2,j_2}(u+r-2) & \cdots & a_{i_2,j_r}(u) \\
          \vdots            &  \vdots            &        &  \vdots     \\
         a_{i_r,j_1}(u+r-1) & a_{i_r,j_2}(u+r-2) & \cdots & a_{i_r,j_r}(u) \\
        \end{bmatrix}
         \notag    \\
=& \sum_{\sigma \in \mathfrak{S}_{r}}
      \sgn{\sigma} 
        a_{i_{\sigma(1)},j_1}(u+r-1) a_{i_{\sigma(2)},j_2}(u+r-2)
                        \cdots a_{i_{\sigma(r)},j_r}(u) 
         \label{e:det(a_+rho-shift)}
\end{align}
if $I=\{ i_1 <i_2 <\cdots <i_r \}$ and 
$J=\{ j_1 <j_2 <\cdots <j_r \}$.

It follows from \eqref{e:def_Xi(u)}, \eqref{e:def_eta},
\eqref{e:anti_commutation_eta} and \eqref{e:Xi_and_eta} that
\begin{align*}
\allowdisplaybreaks
  & \Xi^{(r)}(u+r-1) 
               \\
 =& \, (-)^{r(r-1)/2}
    \sum_{\alpha_1,\alpha_2,\dots,\alpha_r}
    \eta_{\alpha_1}(u+r-1) \eta_{\alpha_2}(u+r-2) \cdots \eta_{\alpha_r}(u)
    e_{-\alpha_1} e_{-\alpha_2} \cdots e_{-\alpha_r}
               \\
 =& \, (-)^{r(r-1)/2}
    \sum_{j_1 <j_2 <\dots <j_r}
       \hspace{3pt}
    \sum_{\sigma \in \mathfrak{S}_{r}}
    \eta_{j_{\sigma(1)}}(u+r-1) \eta_{j_{\sigma(2)}}(u+r-2) 
                                   \cdots \eta_{j_{\sigma(r)}}(u)
               \\
  & \hspace{80pt}
      \times e_{-j_{\sigma(1)}} \cdots e_{-j_{\sigma(r)}}
               \\
 =& \, r! 
    \sum_{j_1 <j_2 <\dots <j_r}
    \eta_{j_1}(u+r-1) \eta_{j_2}(u+r-2) \cdots \eta_{j_r}(u)
                 e_{-j_r} \cdots  e_{-j_1},
               \\
 =& \, r! 
    \sum_{j_1 <j_2 <\dots <j_r}
      \hspace{3pt}
    \sum_{\beta_1, \beta_2, \dots, \beta_r}
       e_{\beta_1} e_{\beta_2} \cdots e_{\beta_r} 
       e_{-j_r} \cdots  e_{-j_1} 
               \\
 & \hspace{80pt}
      \times  a_{\beta_1, j_1}(u+r-1) a_{\beta_2,j_2}(u+r-2) 
                   \cdots a_{\beta_r,j_r}(u) 
              \\
 =& \, r! 
    \sum_{j_1 <j_2 <\dots <j_r}
      \hspace{3pt}
    \sum_{i_1 <i_2 <\dots <i_r}
       e_{i_1} e_{i_2} \cdots e_{i_r} 
       e_{-j_r} \cdots  e_{-j_1} 
               \\
 & \hspace{80pt}
      \times  \sgn{\sigma} \,
             a_{i_{\sigma(1)}, j_1}(u+r-1) a_{i_{\sigma(2)},j_2}(u+r-2) 
                   \cdots a_{i_{\sigma(r)},j_r}(u). 
\end{align*}
This completes the proof.
\qed
\end{proof}

%
%
\begin{lem}
\label{l:Theta2the_power}
For $p,q=0,1,\dots$,
we have
\begin{align*}
 \Theta^p &= 2^p p! \sum_{I \subset [n],|I|=2p}
               e_{I} \Pf{b_I},
                 \\
 \Theta'^q &= 2^q q! \sum_{J \subset [n],|J|=2q}
               e_{-J} \Pf{c_J}.
\end{align*} 
\end{lem}
\begin{proof}
It is a direct calculation to show these formulae, 
as in the proof of Lemma \ref{l:Pfaffian_via_2-form}.
In fact, 
\begin{align*}
\Theta^p 
 &= \sum_{\alpha_1,\beta_1,\dots,\alpha_p,\beta_p}
    e_{\alpha_1} e_{\beta_1} e_{\alpha_2} e_{\beta_2}
                                 \cdots e_{\alpha_p} e_{\beta_p} \,
    b_{\alpha_1,\beta_1} b_{\alpha_2,\beta_2} 
                                 \cdots b_{\alpha_p,\beta_p}
            \notag   \\
 &= \sum_{i_1 <i_2 <\cdots <i_{2p}}
    e_{i_1} e_{i_2} \cdots e_{i_{2p}}
        \hspace{6pt}
    \sum_{\sigma \in \mathfrak{S}_{2p} } 
    \sgn{\sigma} \,
    b_{i_{\sigma(1)},i_{\sigma(2)}} \cdots b_{i_{\sigma(2p-1)},i_{\sigma(2p)}}
                     \\
  &= 2^p p! 
     \sum_{I \subset [n],|I|=2p} e_I \Pf{b_I}.
           \notag 
\end{align*}
The other formula can be shown in a similar way.
\qed
\end{proof}

%
%
\begin{thm}
\label{th:main_thm}
Let 
\(
\bs{X} = \begin{bmatrix}
           a & b \\
           c & -J_n \tp{a} J_n
         \end{bmatrix}
\) 
be the anti-alternating matrix 
defined by \eqref{e:mat_noncommutative}
whose entries are elements of $U(\o_{2n})$,
where $a,b,c$ are given by 
\eqref{e:nc_coloring} and \eqref{e:nc_parameterize}.
Then we have
\begin{equation}
\label{e:nc_msf}
  \Pf{\bs{X}}
   = \sum_{k=0}^{ \lfloor n/2 \rfloor } 
            \sum_{ \substack{I, J \subset [n] \\%
                   |I|=|J| = 2k} }
       \sgn{\bar{I},I} \sgn{\bar{J},J}
        \det \left(  a^{\bar{I}}_{\bar{J}} \,
               +  \mathbbm{1}^{\bar I}_{\bar J} \,\, \bs{\rho}(|\bar J|)
             \right) 
                  \Pf{c_{J}} \Pf{b_{I}},
\end{equation}
where 
$\bs{\rho}(j)$ denotes the diagonal matrix $\diag(j-1,j-2,\dots,1,0)$.
\end{thm}
\begin{proof}
By Propositions \ref{p:trinomial}, \ref{p:Xi2the_power} 
and Lemma \ref{l:Theta2the_power},
we obtain that
\begin{multline}
  \Omega^n = 2^n n! \sum_{ \begin{subarray}{c}
                     p,q,r \geqsl 0 \\
                     p+q+r=n
                    \end{subarray}
                }  \hspace{3pt}
             \sum_{ \begin{subarray}{c}
                     I_1,J_1 \subset [n]  \\
                     |I_1|=|J_1|=r
                    \end{subarray}
                }  \hspace{3pt}
             \sum_{ \begin{subarray}{c}
                     I,J \subset [n]  \\
                     |I|=2p,|J|=2q 
                    \end{subarray}
                }
          e_{I_1} e_{I} e_{-J_1} e_{-J} 
              \notag \\
          \times \det \left(  a^{I_1}_{J_1} \,
               + \mathbbm{1}^{I_1}_{J_1} \diag(u+r-1,u+r-2,\dots,u)
             \right)
          \Pf{c_J} \Pf{b_I}  
\end{multline}
with $u=q-p$.
Since $\Omega^n$ is of top degree,
the terms corresponding to $I_1,I,J_1,J$ in the sum vanish 
unless $I_1 \sqcup I = J_1 \sqcup J = [n]$,
in particular, unless $p=q$.
Thus
\begin{equation*}
\Omega^n=  2^n n! 
        \sum_{k=0}^{\floor{n/2}}
        \sum_{ \begin{subarray}{c}
                 I,J \subset [n]  \\
                 |I|=|J|=2k 
               \end{subarray}
            }
        e_{\bar I} e_{I} e_{-\bar J} e_{-J} 
        \det \left(  a^{\bar I}_{\bar J} \,
               + \mathbbm{1}^{\bar I}_{\bar J} \, \bs{\rho}(n-2k)
             \right)
          \Pf{c_J} \Pf{b_I}.  
\end{equation*}
Now the theorem follows if one compares this 
with the right-hand side of \eqref{e:Omega2the_n}.
\qed
\end{proof}

Using the theorem
one can easily compute 
the eigenvalues of the Pfaffian $\Pf{\bs{X}}$
on the irreducible representaions of $\o_{2n}$.
%
%
%
%
%
\begin{cor}
The eigenvalue of the central element 
$\Pf{\bs{X}} \in ZU(\o_{2n})$  
on the irreducible representaion with highest weight 
$\lambda=\sum_{i=1}^{n} \lambda_i \epsilon_i$
is given by
\( \prod_{i=1}^n (\lambda_i + n - i) \).
\end{cor}
\begin{proof}
Applying $\Pf{\bs{X}}$ to 
the highest weight vector, say $v_{\lambda}$, 
one sees that 
the only term that survives in the sum of \eqref{e:nc_msf} 
is the one corresponding to $I=J=\varnothing$
since each $b_{i,j}$ is a positive root vector
by Remark \ref{r:maximal_parabolic_subalgebra},
and thus by \eqref{e:det(a_+rho-shift)},
one sees that 
\begin{align*}
\label{e:eigenvalue}
\Pf{\bs{X}} v_\lambda 
  &= \det(a + \bs{\rho}(n)) v_\lambda   
       \notag \\
  &= \sum_{\sigma \in \mathfrak{S}_{n}} 
       \sgn{\sigma} a_{\sigma(1),1}(n-1) a_{\sigma(2),2}(n-2)
            \cdots a_{\sigma(n),n}(0) v_\lambda
              \\
  &= (\lambda_1+n-1)(\lambda_2+n-2) \cdots \lambda_n \, v_\lambda,
\end{align*}
since $a_{i,j}$ is also a positive root vector if $i<j$.
\qed
\end{proof}

\subsection*{Acknowledgements}
I would like to thank Professor T\^{o}ru Umeda 
for drawing my attention to \cite{IU01}.

%
%
\appendix

\section{Proof of Commutative Minor Summation Formula via Exterior Algebra}

In this Appendix,
we give another proof of
the commutative minor summation formula \eqref{e:msf} 
in Theorem \ref{th:msf},
using the exterior calculus.

Let 
$X=(x_{i,j})_{i,j} \in \o_{2n}$
be an anti-alternating matrix 
with commutative entries:
\begin{equation*}
X
=\left[
 \begin{array}{cccc|ccc}
  x_{1,1} & x_{1,2} & \cdots & x_{1,q} & x_{1,-p} & \cdots & 0         
    \\
  \vdots  & \vdots  &        & \vdots  & \vdots   & \adots & \vdots          
    \\  
  x_{p,1} & x_{p,2} & \cdots & x_{p,q} & 0        & \cdots & x_{p,-1} 
    \\[3pt]  \hline
  x_{-q,1}& x_{-q,2}& \cdots & 0       &x_{-q,-p} & \cdots & x_{-q,-1}
    \\[2pt]
  \vdots  & \vdots  & \adots & \vdots  & \vdots   &        & \vdots
    \\[2pt]
  x_{-2,1}& 0       & \cdots & x_{-2,q}&x_{-2,-p} & \cdots & x_{-2,-1}
    \\[3pt]
  0       & x_{-1,2}& \cdots & x_{-1,q}&x_{-1,-p} & \cdots & x_{-1,-1}
 \end{array}
\right], 
\end{equation*}
where $x_{-j,-i}=-x_{i,j}$ for all $i,j$.
Define 2-forms $\Omega$ by
\[
 \Omega
  := \sum_{ \begin{subarray}{c}
               i \,\in \,[p] \cup [-q]  \\
               j \,\in \,[q] \cup [-p] 
            \end{subarray}
        }
     e_i e_{-j} x_{i,j}
\]
for $X$.
By the same argument as in the proof of 
Lemma \ref{l:Pfaffian_via_2-form} 
one can show that
\begin{equation*}
 \Omega^n 
 = e_1 \cdots e_p e_{-q} \cdots e_{-1} 2^n n! \Pf{X}. 
\end{equation*}
Note that by our convention
$-q$ stands for $2n+1-q=p+1$ and so on.
Coloring and parameterize $X$ as in
\eqref{e:coloring} and \eqref{e:parameterize},
define 2-forms 
$\Xi,\Theta,\Theta'$ by
\begin{equation*}
\Xi = \sum_{i \in [p], j \in [q]}    e_i e_{-j} \, a_{i,j},  
     \quad  
\Theta = \sum_{i,j \in [p]}  e_i e_{j} \, b_{i,j},  
     \quad
\Theta' = \sum_{i,j \in [q]} e_{-j} e_{-i} \, c_{i,j}.  
\end{equation*}
It is clear that
\[
 \Omega = \Theta' + 2 \Xi + \Theta.
\]
Furthermore,
the trinomial expansion formula 
in Proposition \ref{p:trinomial} holds
without parameter-shift:
\begin{equation}  
\label{e:commutative_trinomial}
 \Omega^m = \sum_{ \begin{subarray}{c}
                   h,s,t \geqsl 0 \\
                   h+s+t=m
                 \end{subarray}
               }
           \frac{m!}{h! s! t!} 2^{h} \, 
             {\Xi}^{h} \Theta^{s} \Theta'^{t}
\end{equation}
for $m=0,1,\dots$.
The same calculation as in 
Proposition \ref{p:Xi2the_power} and Lemma \ref{l:Theta2the_power}
yields
\begin{align*}
\Xi^{h}
 &= h! \sum_{ \begin{subarray}{c}
             I \subset [p], J \subset [q] \\
             |I|=|J|=h
            \end{subarray}
        }
     e_{I} e_{-J} \det ( a^{I}_{J} ),
                 \\ 
\Theta^{s} 
 &= 2^s s! \sum_{I \subset [p],|I|=2s}
               e_{I} \Pf{b_I},
                 \\
{\Theta'}^t 
 &= 2^t t! \sum_{J \subset [q],|J|=2t}
               e_{-J} \Pf{c_J}
\end{align*}
for $h,s,t=0,1,\dots$.
Substituting these into \eqref{e:commutative_trinomial}
with $m=n$,
we see that
\begin{equation}   
\label{e:commutative_Omega2then_n}
  \Omega^n = 2^n n! \sum_{ \begin{subarray}{c}
                     h,s,t \geqsl 0 \\
                     h+s+t=n
                    \end{subarray}
                }  \hspace{3pt}
             \sum_{ \begin{subarray}{c}
                     I_1,I \subset [p], J_1,J \subset [q]  \\
                     |I_1|=|J_1|=h    \\
                     |I|=2s, |J|=2t 
                    \end{subarray}
                }  
          e_{I_1} e_{I} e_{-J_1} e_{-J} 
              \det (a^{I_1}_{J_1})  \Pf{b_I}  \Pf{c_J}.
\end{equation}
Since $\Omega^n$ is of top degree,
the terms that survive in the sum \eqref{e:commutative_Omega2then_n}
are those corresponding to $I_1,I,J_1,J$ 
satisfying $I_1 \sqcup I=[p]$ and $J_1 \sqcup J=[q]$.
In particular,
$h+2s=p$ and $h+2t=q$, 
where $|I_1|=|J_1|=h, |I|=2s, |J|=2t$.
Now setting $h=2k+\epsilon$ 
with $\epsilon$ the parity of $p$ (=the parity of $q$),
we obtain the formula.


\bibliographystyle{amsalpha}
\bibliography{rep}

\nocite{MN99}
\nocite{HU91}

\end{document}